%% file: relhyp.tex
\documentclass[12pt]{article}

\usepackage{amsmath,amssymb,amsthm,amsfonts,graphicx,color}
\usepackage{epsfig}

\newtheorem{thm}{Theorem}[section]

\newtheorem{cor}[thm]{Corollary}
\newtheorem{prop}[thm]{Proposition}

\theoremstyle{remark} 

\newtheorem{remark}[thm]{Remark}

\newtheorem{definition}[thm]{Definition}

\def\H{{\mathbb H}}

\def\g{{\Gamma}}

\title{A combination theorem for relatively hyperbolic groups}
\date{\today}
\author{Emina Alibegovi\'{c}}

\begin{document}
\renewcommand{\thefootnote}{\null}
\maketitle
\setcounter{footnote}{0}
\renewcommand{\thefootnote}{\arabic{footnote}}

\begin{abstract}
In this paper we give new requirements that a tree of $\delta$-hyperbolic
spaces has to satisfy in order to be $\delta$-hyperbolic itself. As an
application, we give a simple proof that limit groups are relatively
hyperbolic. 
\end{abstract}

\section{Introduction}

In his work on Diophantine equations over free groups Z. Sela
introduced limit groups. He showed that this class of groups
coincides with the class of $\omega$-residually free groups. Since then
few more characterizations of limit groups have been given, as well as
structure theorem for limit groups. 

This work introduced a lot of interesting questions one might ask
about limit groups.  
We were interested in describing the set of homomorphisms from an
arbitrary f.g. group $G$ into a limit group $L$, $Hom(G,L)$. A key tool
in studying $Hom(G,L)$ is a $\delta$-hyperbolic space on
which the given limit group $L$ acts freely, by isometries. We
construct such a space in Section \ref{3}.
That the space we constructed is $\delta$-hyperbolic follows from
Theorem \ref{combination} which  gives conditions a tree of
hyperbolic spaces has to satisfy in order to be a hyperbolic space.
The proof of this theorem is given in Section \ref{2} and is an
adaptation of the proof of Bestvina-Feighn Combination Theorem to a
different setting. 

Existence of such spaces for limit groups gives an answer to the
question whether limit groups are hyperbolic relative to their maximal
noncyclic abelian subgroups. This question was answered
affirmatively by F. Dahmani, \cite{combconv}, who proved a combination
theorem for geometrically finite convergence groups using different methods. 

{\bf Acknowledgements} I would like to thank Mladen Bestvina for his
support. I will always be grateful for the time and knowledge he generously
shared with me.

\section{A combination theorem}\label{2}

Let $X$ be a connected finite cell complex which is a graph of
spaces. There is a map $p:X \to \g$ onto a finite graph $\g$. Let
$X_e$ denote the preimage of a midpoint of an edge $e$ in $\g$, and
let $X_v$ denote the preimage of a component of $\g \backslash \{
\text{midpoints of all edges}\}$ that contains the
vertex $v$. We require that $X_e$ and $X_v$ are connected and that
their inclusions into $X$ induce inclusions on fundamental
groups. There is an induced map $\tilde{p}: \tilde{X} \to T$ from the
universal cover of $X$ onto a $\pi_1(X)$-tree $T$ so that $T/\pi_1(X)$
is isomorphic to $\g$. We say that $X$ is a {\it graph of negatively
  curved spaces} if every vertex space $X_v$ is negatively curved. As
a reminder a cell complex $X$ is said to be negatively curved if there
exists a constant $A=A(X)$ so that each inessential circuit bounds a
disk of combinatorial area which is bounded above by $A$ times the
combinatorial length of the circuit (we assign length 1 to each
edge). 

We know that in $\delta$-hyperbolic spaces geodesic
triangles are $\delta$-thin. We have a similar fact for polygons, in
fact for quasigeodesic polygons. The following
proposition can be found in \cite{gromov} and \cite{combi}. 

\begin{prop}\label{qipolygons}
Let $Z$ be a $\delta$-hyperbolic space and let $\tau \geq 1$ be a
constant. There is a function $B(x)=O(\log x)$ and a linear function
$C(x)$ each depending only on $Z$ and $\tau$ with the following
property. If $\Delta: D^2 \to Z$ is a disk with boundary a $k$-sided
$\tau$-quasigeodesic polygon, then there is a finite tree $S$ and a
map $r:D^2 \to S$ such that:
\begin{enumerate}
\item the number of valence one vertices of $S$ is $k$,
\item for $a$ and $b$ in $S^1$, $d_Z(\Delta(a),\Delta(b))\leq
  d_S(r(a),r(b)) + B(k)$,
\item $r^{-1}(s)$ is a properly embedded finite tree in $D^2$ for $s
  \in S$, 
\item if $E$ is an edge of $S$, then $r$ restricted to
  $r^{-1}(\text{Interior}(E))$ is an $I$-bundle,
\item for $a_1,b_1$ (respectively $a_2,b_2$) in the same side of the
  polygon and satisfying $r(a_1)=r(a_2) \in E$ and $r(b_1)=r(b_2) \in
  E$, we have
\[\ell(\Delta(\, \text{the circular arc}\ a_1b_1 \ \text{in the edge
   of the polygon}\ ) \leq\]
\[C(\ell(\Delta(\, \text{the circular arc} \ a_2b_2 \ \text{in the edge
    of the polygon}\ ))).\]
\end{enumerate}
\end{prop}
\noindent
Such a map $r$ is called a resolution of the quasigeodesic polygon. A
singular fiber of the resolution is a fiber which is not isomorphic to
$I$. 

{\it Partial qi-embedded condition} We will say that a graph of spaces
$X$ satisfies the partial quasiisometrically embedded condition if
every edge space $\widetilde{X_e}$ is quasiisometrically embedded in
at least one of the vertex spaces $\widetilde{X_v}$ and
$\widetilde{X_w}$, where $v$ and $w$ are endpoints of the edge $e$ in
$\g$. We further ask that all qi-constants be equal. 

{\it Qi-consistency condition} A graph of groups $X$ satisfies the
qi-consistency condition if the following holds: If one of the
edge spaces adjacent to a vertex space $\widetilde{X_v}$ qi-embeds
into it, then the same is true for all adjacent edge spaces. We will
call such vertex spaces {\it good}. 

{\it Compact intersection condition} A graph of spaces $X$ satisfies
the compact intersection condition if whenever the edge spaces
$\widetilde{X_e}$ and $\widetilde{X_f}$ qi-embed into the same vertex
space $\widetilde{X_v}$ then the intersection of any of their
Hausdorff neighborhoods in $\widetilde{X_v}$ is a compact set.

\begin{remark}\label{L} The compact intersection condition gives us
  the following. Suppose that $\widetilde{X_e}$ and $\widetilde{X_f}$
  qi-embed into the same vertex space $\widetilde{X_v}$. If we fix a
  constant $k < \infty$, then  
\[L=\max\{\ell(S_1) : \exists S_2 \subset \widetilde{X_f}\, \text{so
  that}\, d_H(S_1,S_2) \leq k \} \]
is finite, where $S_1, S_2$ are quasigeodesics in $\widetilde{X_e}$
and $\widetilde{X_f}$, respectively. Note that $L$ depends on the
qi-constants for $S_1$ and $S_2$. 
\end{remark}
%\begin{thm}\label{combination}
%Let $X$ be a tree of $\delta$-hyperbolic spaces some of that we
%designate as 'nice'
%vertex spaces so that the following conditions are satisfied:
%\begin{itemize}
%\item Every edge space is incident with at least one 'nice' vertex
%  space. Further, it is (quasi)isometrically embedded into all the
%  adjacent 'nice' vertex spaces.
%\item For every two edge spaces that are embedded into the same
%  'nice' vertex space the intersection of any of their Hausdorff
%  neighborhoods is a compact set.
%\end{itemize}
% $X$ is then a hyperbolic space. 
%
%\end{thm}

\begin{thm}\label{combination}
If a tree of negatively curved spaces $X$ satisfies the partial
qi-embedded, the qi-consistency and the compact intersection conditions, then
$X$ is negatively curved.  

\end{thm}

We would like to show that $X$ satisfies subquadratic isoperimetric
inequality, which will then imply that $X$ is a hyperbolic space
(\cite{bowditchsub}, \cite{gromov}). We will use the techniques
employed by Bestvina and Feighn in the proof of their Combination
theorem (\cite{combi}).

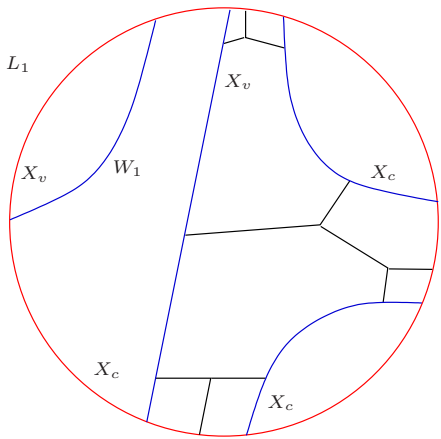
\begin{figure}[!ht]
\centerline{\input{walls1.pstex_t}}
\caption{Decomposition by walls of a disk $\Delta$ into polygons}\label{walls}
\end{figure}

Let $\gamma:S^1 \to X$ be a circuit that is transverse to and has
nonempty intersection with $\cup\{X_e: e \ \text{edge of} \ T\}$. We
may also assume that $\gamma$ is contained in the 1-skeleton of $X$
(\cite{jenpepe}). Following \cite{combi} we talk about {\it good
disks}. There is a disk $\Delta:D^2 \to X$ with boundary $\gamma$. The set
$\mathcal{W}=\Delta^{-1}(\cup\{X_e: e\ \text{edge of} \ T\})$ divides
$D^2$ into regions that are mapped into negatively curved vertex
spaces, see Figure \ref{walls}. Elements of $\mathcal{W}$ are called
{\it walls}. We may assume that $\Delta$ has the following properties:
\begin{enumerate}

\item[(1)] The set $\mathcal{W}$ consists of properly embedded arcs in
  $D^2$. 
\item[(2)] The length of $\Delta(\cup \mathcal{W})$ in $X$ is
  minimal over all disks satisfying (1).
\item[(3)] The closures of the components of
  $\Delta(D^2\backslash (\cup \mathcal{W}))$ have areas bounded by $A$
  times the length of their boundaries, where $A$ is a constant. 
\item[(4)] Define $\mathcal{L}$ to be the set of closures of the components
  of $S^1\backslash (S^1 \cap \cup \mathcal{W})$. We may assume that
  $\gamma$ restricted to each element of $\mathcal{L}$ is a geodesic in the
  appropriate $X_v$. We view $\gamma$ as a polygon whose
  sides are elements of $\mathcal{L}$. Hence the number of sides of
  $\gamma$ can be no more than $\ell(\gamma)$ (the length of each side
  is at least 1).
\end{enumerate}
\noindent
A disk is good if it satisfies (1)-(4). 

\noindent

Our goal is to bound the area of $\Delta$ by a subquadratic function
of the length of its boundary. For a good disk $\Delta$, by property
(3), we have
$$Area(\Delta) \leq A\,(2\ell(\Delta(\cup
\mathcal{W}))+\ell(\Delta(\cup \mathcal{L}))) \leq A\,(2\ell(\Delta(\cup
\mathcal{W}))+\ell(\gamma))).$$
Therefore we need to bound $\ell(\Delta(\cup\mathcal{W}))$ in terms of
$\ell(\gamma)$. 

Let us denote by $\mathcal{P}$ the set of closures of the components
of $D^2\backslash \cup \mathcal{W}$, and let $P$ be an element of
$\mathcal{P}$. If $\Delta(P)$ is contained in a good vertex space
$X_v$, then the map $\Delta$ restricted to $\partial P$ is a
$\tau$-quasigeodesic polygon in $X_v$. Note that each wall $W \in
\mathcal{W}$ is a side of at least one polygon $P\in \mathcal{P}$ for
which $\Delta(P)$ is contained in a good vertex space, and so we only
need to consider such polygons. In what follows in order to avoid a
cumbersome notation we will write $\ell(W)$ and $P$ when we really
mean $\ell(\Delta(W))$ and $\Delta(P)$.

{\it Claim 1:} 
\[\sum \{\ell(W): W\subset \, \text{a bigon} \, P \subset X_v\} \leq
\tau \ell(\gamma). \]

{\it Proof of Claim} If $P$ is a bigon whose one side is a wall $W$,
then the other side $s$ of that bigon is an element of
$\mathcal{L}$. According to (4) in the definition of good disks the
images of the elements of $\mathcal{L}$ under $\Delta$ are geodesics
in appropriate vertex spaces and hence $\ell(W) \leq \tau \ell(s)$. The
claim follows.

{\it Claim 2:} 
\[\sum \{\ell(W): W\subset \,\text{an}\,m \,\text{-gon}\,P, m\geq
4\}=O(\ell(\gamma)\log(\ell(\gamma))).\] 

{\it Proof of Claim} A polygon $P$ is a $\tau$-quasigeodesic polygon
and can be resolved using Lemma \ref{qipolygons}.
We call a point $w \in W$ a singular point if it lies in a singular
fiber of the resolution of the polygon $P$. These points will
decompose $W$ into the union of closed segments $V$.  
%$$W=\cup\{V:\, \text{endpoints of}\, V \, \text{are singular points or
%  they lie on}\, \Delta(S^1) \}.$$
%\noindent
Let $\mathcal{V}(W)$ denote the set of all such $V$'s, and let
$\mathcal{V}=\cup \{\mathcal{V}(W): W \in \mathcal{W}\}$. Since 
$$\sum \{\ell(W): W \subset\,\text{an}\,m
  \,\text{-gon}\,P, m\geq 4\}=\sum \{\ell(V): V \in \mathcal{V}\}$$
\noindent
we need to bound $\ell(V)$, for all $V \in \mathcal{V}$.

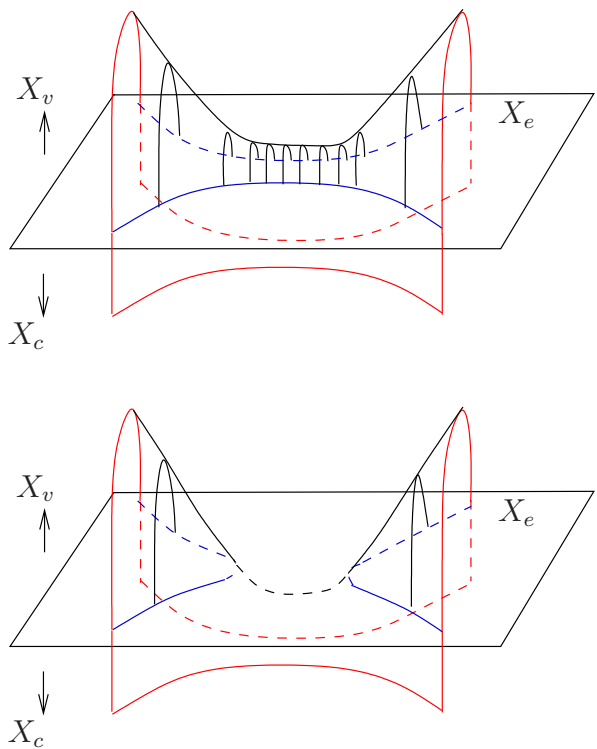
\begin{figure}[!ht]
\centerline{\input{surgery.pstex_t}}
\caption{Surgery: the lengths of the central 'parallel' curves are
  much shorter than the lengths of walls connecting them. We change a
  disk by pushing this tunnel down into $X_c$ and get a new disk: a
  good disk. For the better visibility we drew $X_e$ as two
  dimensional.}\label{surg}
\end{figure}

Note that singular fibers decompose the polygon $P$ into the union of
quadrilaterals and triangles. The case of triangles is relatively easy
to handle.  One of the sides of this triangle is some $V \in
\mathcal{V}$ and another, call it $S$, is contained in $S^1$.
Considering how the resolution was formed, we have that $\ell(V)\leq
C(\ell(S))$, where the function $C$ is a linear function from
Proposition \ref{qipolygons} (5). Hence, 
\[\sum \{\ell(V): V \in \mathcal{V}\}\leq D(\ell(\gamma))+\sum
\{\ell(V): V\, \text{side of quadrilateral in}\, P\}\] where $D$ is a
linear function.  Let us consider the case of a quasigeodesic
quadrilateral $Q$ with sides $V_1$ and $V_2$ (contained in
$X_e$ and $X_f$, respectively) that are joined
by singular fibers of the resolution $r$ of the polygon $P$. We may
assume that $V_1$ is shorter of the two.  According to Proposition
\ref{qipolygons} the distance between the images under $\Delta$ of the
two endpoints of a fiber is at most $B(\, \text{the number of sides
of}\, P) \leq B(\ell(\gamma))$.  Since $Q$ is contained in
$\delta$-hyperbolic space, it is not hard to show that
as long as $\ell(V_1) > 2\tau B(\ell(\gamma))+d$, where $d=d(\tau,
\epsilon, R)$, we can find subsegments $S_1 \subset V_1$ and $S_2
\subset V_2$ such that $d_H(S_1,S_2)\leq 12 \delta + 4R(\tau
+1)+2\epsilon$.  If $\ell(S_1)>L$, ($L$ from Remark \ref{L}, where the
constant $k$ is $12 \delta + 4R(\tau +1)+2\epsilon$), then both $S_1$
and $S_2$ belong to $X_e$. By doing a surgery, see Figure
\ref{surg}, we can shorten the length of the walls in our disk, which
contradicts the assumption that $\Delta$ was a good disk.

Thus, if
\[\ell(V) > 2\tau B(\ell(\gamma))+d +L, \ \text{for some}\ V
\in \mathcal{V},\] 
$\Delta$ is not a good disk.
\noindent
To finish the proof we need to know what the cardinality of
$\mathcal{V}$ is. The number of $V \in \mathcal{V}$ is proportional to the
number of singular fibers inside $\Delta$. On the other hand, there can
be no more singular fibers than there are triangles in the
triangulation of the polygon $\gamma$. Hence,
cardinality$(\mathcal{V})=O(\ell(\gamma))$,
%. Finally, 
%$$\ell(\Delta(\cup \mathcal(W)))=O(\ell(\gamma)\log(\ell(\gamma))),$$
and our claim follows. 

\begin{proof}[Proof of Theorem \ref{combination}] 
If $\Delta$ is a good disk we already noticed that 
\begin{eqnarray*}
\text{Area}(\Delta)& \leq & A\,(2\ell(\Delta(\cup
\mathcal{W}))+\ell(\Delta(\cup \mathcal{L})))\\
& \leq & A\,(2\ell(\Delta(\cup \mathcal{W}))+\ell(\gamma)))\\
%& = & 2A\sum \{\ell(W): W\subset \, \text{bigon} \, P\in \mathcal{P}\}+ 
%2A\sum \{\ell(W): W\subset m \,\text{-gon}\,P, m\geq
%3\}+A\ell(\gamma)\\
&\leq& 2A (\tau\ell(\gamma)+O(\ell(\gamma)\log(\ell(\gamma))))+A\ell(\gamma)\\
&=& O(\ell(\gamma)\log(\ell(\gamma)))
\end{eqnarray*}
The last inequality follows from Claims 1 and 2. Therefore,
$X$ satisfies subquadratic isoperimetric inequality, and our
proof is finshed.
\end{proof}

\section{Application to limit groups}\label{3}

The goal of this section is to produce a $\delta$-hyperbolic space on
which a given limit group acts freely, by isometries. We will in fact
consider a slightly larger class of groups, $\mathcal{C}$. We describe
elements of $\mathcal{C}$ ($\mathcal C$-groups, for short)
inductively. 

\begin{definition}
A torsion-free, f.g. group $G$ is a depth 0 $\mathcal
C$-group if it is either an f.g. free group, or an f.g. free abelian
group or the fundamental group of a closed hyperbolic surface. A torsion-free
f.g. group $G$ is a $\mathcal C$-group of depth $\leq n$ if it has a graph of
groups decomposition with three types of vertices: abelian, surface or
depth $\leq (n-1)$, cyclic edge stabilizers and the following holds:

\begin{itemize}
\item Every edge is adjacent to at most one abelian vertex $v$.
  Further, $G_v$, the stabilizer of $v$, is a maximal abelian subgroup
  of $G$.
\item Each surface vertex group is the fundamental group of a surface with
  boundary, and to each boundary component corresponds an edge of this
  decomposition. Each edge group is conjugate to a boundary
  component. 
\item The stabilizer of a depth $\leq (n-1)$ vertex $v$, $G_v$, is
  $\mathcal C$-group of depth $\leq (n-1)$. The images in $G_v$ of
  incident edge groups are distinct maximal abelian subgroups of $G_v$
  (i.e., cyclic subgroups generated by distinct, primitive, hyperbolic
  elements of $G_v$).
\end{itemize}

We say that the depth of a $C$-group $G$ is the smallest $n$ for which
$G$ is of depth $\leq n$. 
\end{definition}

We will get a hyperbolic space on which a $\mathcal C$-group $G$ acts
freely by induction on its depth.  If $G$ is a depth 0 $\mathcal
C$-group, we take a tree, a horoball or $\H^2$.
% We remark that if $G=\Z^n$
%then the space we want is a horoball in $\H^{n+1}$. We will take a
%horoball in an upper-half space model, so that $G$ acts by Euclidean
%translations on hyperplanes parallel to $\R^n$.

Let $G$ be a depth $n$ $\mathcal C$-group. For the vertex
groups in the decomposition of $G$ as above we have the desired spaces by
induction. Let $X/G$ be a graph of spaces corresponding to this
splitting, and $X$ its universal cover.  Let $T_{G}$ be the
corresponding graph of groups, and let $T$ be a tree so that
$T/G=T_G$. Our goal is to show that $X$ satisfies the hypotheses of
Theorem \ref{combination} and consequently that $X$ is a hyperbolic space.

The requirement we imposed on the splitting of $G$ guarantees that
$X/G$ satisfies the partial qi-embedded condition and the qi-consistency
condition.  Namely, the generators of all edge groups adjacent to a
nonabelian vertex group are identified with hyperbolic elements of
that vertex group, hence the corresponding edge spaces in $X$ are
glued along quasigeodesics in the relevant vertex space. Just as a
note, no edge space qi-embeds into a vertex space which is a horoball.

{\it Claim:} $X/G$ satisfies the compact intersection condition. 

{\it Proof of Claim:} Suppose edge spaces $\widetilde{X_e}$ and
$\widetilde{X_f}$ qi-embed into the vertex space $\widetilde{X_v}$. We
noticed that they are glued along quasigeodesics, say $c_1$ and $c_2$
respectively. There are elements $g_1$ and $g_2$ of $\pi_1{X_v}$ that
act as translations along $c_1$ and $c_2$, respectively.  If the
Hausdorff distance between quasigeodesics $c_1$ and $c_2$ is bounded,
we conclude that $g_1$ and $g_2$ fix the same two points in $\partial
\widetilde{X_v}$. Hence, they are both contained in a unique elementary group that
is virtually cyclic. Since we have no torsion elements, this
elementary group is cyclic, contradicting the choice of splitting for
$G$.
\qed\\

The definition of relatively hyperbolic groups appears in
many forms. We use the one given by M. Gromov in \cite{gromov}, 8.6.
 
Let $X$ be a complete hyperbolic locally compact geodesic space with a
discrete free isometric action of a group $\g$ such that the quotient
$V=X/\g$ is quasiisometric to the union of $k$ copies of $[0,\infty)$
joined at 0. Lift $k$ rays that correspond to $\partial V$ to rays
$r_i:[0,\infty) \to X, \ i=1,\ldots,k$. Let $h_i$ be the horofunction
corresponding to $r_i$ and let $r_i(\infty)$ be the limit point of
$r_i$. Denote by $\g_i<\g$ the stabilizer of $r_i(\infty)$ and assume
that it preserves $h_i$. Denote by $B_i(\rho)$ the horoballs
$h_i^{-1}(-\infty,\rho)\subset X$ and assume that for sufficiently
small $\rho$ the intersections $\gamma B_i(\rho) \cap B_j(\rho)$ is empty
unless $i=j$ and $\gamma \in \g_i$. Let 
\[\g B(\rho)=\bigcup_{i, \gamma}\gamma B_i(\rho), \] 
$i=1, \ldots, k, \ \gamma \in \g$. Let $X(\rho)=X\backslash\g B(\rho)$, and
assume that $X(\rho)/\g$ is compact for all $\rho \in (-\infty, \infty)$. 

\begin{definition}\label{rhgroups}
We say that a group $\g$ is hyperbolic relative to the subgroups
$\g_1, \ldots, \g_k$ if $\g$ admits an action on some $X$ as above,
and where $\g_i$ are the stabilizers of $h_i$. 
\end{definition}

After inspection of the action of $G\in \mathcal C$ on the
$\delta$-hyperbolic space $X$ that we constructed above, we see that
all the requirements of Definition \ref{rhgroups} are
satisfied. Hence, we have proved: 

\begin{thm}\label{rh}
Groups in $\mathcal C$ are hyperbolic relative to the collection of the
conjugacy classes of their maximal noncyclic abelian subgroups.
\end{thm}

Since limit groups belong to the class $\mathcal C$, see \cite{zlil1}
Theorem 3.2. and Theorem 4.1., the consequence of this theorem is the 
following corollary: 

\begin{cor}
Limit group $L$ is hyperbolic relative to the collection of
representatives of conjugacy classes of its maximal noncyclic abelian
subgroups. 
\end{cor}

This fact was also noted in the work of F. Dahmani, \cite{combconv}. 

%We already mentioned that there are multiple definitions of relative
%hyperbolicity in the literature. 
Several nice properties follow from the relative
hyperbolicity. I. Bumagin  showed that the conjugacy problem is solvable
for a group $G$ which is hyperbolic relative to a subgroup $H$ with
solvable conjugacy problem; hence limit groups have solvable conjugacy
problem (\cite{bumagina}). D. Y. Rebecchi showed that a group $G$
hyperbolic relative to a biautomatic subgroup $H$ is itself
biautomatic, \cite{rebbechi}. We therefore conclude that limit groups
are biautomatic.

\bibliographystyle{siam}
\bibliography{../global/all}
\par

\end{document}

%% file: walls1.pstex_t
\begin{picture}(0,0)%
\epsfig{file=walls1.pstex}%
\end{picture}%
\setlength{\unitlength}{3947sp}%
\begingroup\makeatletter\ifx\SetFigFont\undefined%
\gdef\SetFigFont#1#2#3#4#5{%
  \reset@font\fontsize{#1}{#2pt}%
  \fontfamily{#3}\fontseries{#4}\fontshape{#5}%
  \selectfont}%
\fi\endgroup%
\begin{picture}(2714,2706)(2234,-4305)
\put(2236,-1984){\makebox(0,0)[lb]{\smash{\SetFigFont{7}{8.4}{\rmdefault}{\mddefault}{\itdefault}\(L_1\)}}}
\put(2909,-2639){\makebox(0,0)[lb]{\smash{\SetFigFont{7}{8.4}{\rmdefault}{\mddefault}{\itdefault}\(W_1\)}}}
\put(2329,-2676){\makebox(0,0)[lb]{\smash{\SetFigFont{7}{8.4}{\rmdefault}{\mddefault}{\itdefault}\(X_v\)}}}
\put(2788,-3910){\makebox(0,0)[lb]{\smash{\SetFigFont{7}{8.4}{\rmdefault}{\mddefault}{\itdefault}\(X_c\)}}}
\put(3609,-2096){\makebox(0,0)[lb]{\smash{\SetFigFont{7}{8.4}{\rmdefault}{\mddefault}{\itdefault}\(X_v\)}}}
\put(4525,-2667){\makebox(0,0)[lb]{\smash{\SetFigFont{7}{8.4}{\rmdefault}{\mddefault}{\itdefault}\(X_c\)}}}
\put(3881,-4115){\makebox(0,0)[lb]{\smash{\SetFigFont{7}{8.4}{\rmdefault}{\mddefault}{\itdefault}\(X_c\)}}}
\end{picture}

%% file: surgery.pstex_t
\begin{picture}(0,0)%
\epsfig{file=surgery.pstex}%
\end{picture}%
\setlength{\unitlength}{3947sp}%
\begingroup\makeatletter\ifx\SetFigFont\undefined%
\gdef\SetFigFont#1#2#3#4#5{%
  \reset@font\fontsize{#1}{#2pt}%
  \fontfamily{#3}\fontseries{#4}\fontshape{#5}%
  \selectfont}%
\fi\endgroup%
\begin{picture}(3691,4612)(49,-5721)
\put(116,-1676){\makebox(0,0)[lb]{\smash{\SetFigFont{12}{14.4}{\rmdefault}{\mddefault}{\itdefault}\(X_v\)}}}
\put( 61,-3221){\makebox(0,0)[lb]{\smash{\SetFigFont{12}{14.4}{\rmdefault}{\mddefault}{\itdefault}\(X_c\)}}}
\put(3141,-1831){\makebox(0,0)[lb]{\smash{\SetFigFont{12}{14.4}{\rmdefault}{\mddefault}{\itdefault}\(X_e\)}}}
\put(116,-4176){\makebox(0,0)[lb]{\smash{\SetFigFont{12}{14.4}{\rmdefault}{\mddefault}{\itdefault}\(X_v\)}}}
\put(3141,-4331){\makebox(0,0)[lb]{\smash{\SetFigFont{12}{14.4}{\rmdefault}{\mddefault}{\itdefault}\(X_e\)}}}
\put( 61,-5721){\makebox(0,0)[lb]{\smash{\SetFigFont{12}{14.4}{\rmdefault}{\mddefault}{\itdefault}\(X_c\)}}}
\end{picture}